\documentclass[10pt]{article}
\usepackage{amsfonts}
\usepackage{amssymb}
\usepackage{epsfig}
\usepackage{lscape}
\usepackage{comment}
\input{epsf.sty}
\input xypic
\xyoption{all}
\LaTeXdiagrams
\xyoption{2cell}
\UseAllTwocells

\newcommand{\f}[2]{\frac{\displaystyle #1}{\displaystyle #2}}
\def\sq{\sqrt}
\def\sq2{\sqrt{2}}
\def\sq12{\sq{12}}
\def\dsq2{\f{1}{\sqrt{2}}}
\def\be{\begin{equation}}
\def\ee{\end{equation}}
\def\q={\quad = \quad}
\def\lra{\longrightarrow}
\def\Lra{\Longrightarrow}

\def\C{{\mathbb C}}

\def\Set{{\bf Set}}

\newcounter{examnum}[section]
\newcounter{remarnum}[section]
\begin{document}
\title{Mathematics via Symmetry}
\author{Noson S. Yanofsky\footnote{Department of Computer and Information Science,
Brooklyn College, CUNY, Brooklyn, N.Y. 11210 and the Computer Science
Department of the Graduate Center, CUNY, New York, N.Y. 10016. e-mail:
noson@sci.brooklyn.cuny.edu}\qquad Mark Zelcer\footnote{e-mail: mzelcer1@gc.cuny.edu}}
\maketitle
\begin{abstract} 
\noindent We state the defining characteristic of mathematics as a type of symmetry where one can change the connotation of a mathematical statement in a certain way when the statement's truth value remains the same. This view of mathematics as satisfying such symmetry places mathematics as comparable with modern views of physics and science where, over the past century, symmetry also plays a defining role. We explore the very nature of mathematics and its relationship with natural science from this perspective. This point of view helps clarify some standard problems in the philosophy of mathematics. 
\end{abstract}

\section{Introduction}

For\footnote{We wish to thank our friend and mentor Distinguished Professor Rohit Parikh for helpful conversations and for much warm encouragement. We have gained much from his open minded scholarship and mathematical and philosophical sophistication. The first author would also like to thank Jim Cox and Dayton Clark for many stimulating conversations on some of these topics. Thanks also to Jody Azzouni, Sorin Bangu, Andr\'{e} Lebel, Jim Lambek, Guisseppe Longo, Jolly Mathen, Alan Stearns, Andrei Rodin, Mark Steiner, Robert Seely, and Gavriel Yarmish who were extremely helpful commenting on earlier drafts.} a number of years philosophers of science have been making use of symmetry and arguments that emerge from considerations of symmetry to understand how science works and to explicate fundamental concepts in physical theories. Most related philosophical questions deal with specific symmetries, objectivity, interpreting limits on physical theories, classification, and laws of nature (see \cite{Brading2003} for examples of each). In this essay we argue that symmetry also lies at the heart of the nature of mathematics. 

The recent history of the philosophy of mathematics is largely focused on grasping and defining the nature and essence of mathematics and its objects. Attempts to do this include explicating versions of: mathematics is just logic, mathematics is just structure, mathematics is a meaningless game, mathematics is a creation of the mind, mathematics is a useful fiction, etc. 

We believe that a successful account of the nature of mathematics must accomplish a number of goals. First it must square with practice. It cannot ignore the work of any large segment of the mathematical community, past or present. Nor should a credible account of mathematics divide mathematics along lines that make little intuitive sense to mathematicians, such as the philosophical divisions between literal and figurative mathematics, or mathematics to which we have more or less epistemic access, etc. 

Second, a successful account must avoid the main problems already within traditional philosophy of mathematics. It must, for example, account for the usefulness of mathematics in ways that formalists are sometimes accused of being unable to do. It must describe the  nature of mathematical objects in the way that numbers-as-$x$ theories cannot do well (e.g. numbers as sets, mental objects, forms, etc). It should also not force us to countenance new levels of abstracta or novel objects, as some realisms require, merely to satisfy philosophical scruples. 

It is also advantageous if an approach is consistent and continuous with {\em scientific} practice. We would do well to join all of our scientific postulates with our mathematical statements under a uniform semantic, metaphysical, and/or methodological schema while simultaneously addressing traditional questions that emerge from the apparent disconnect between the ontology and methodology of science and mathematics. 

Here we offer a definition of mathematics as a subject matter whose nature completely emerges from symmetry considerations. This conception will satisfy most of the desiderata just mentioned. Just as symmetry considerations have helped clarify many issues in the philosophy of science, so too our work sheds light on familiar problems in the philosophy of mathematics. 

The argument for our definition of mathematics reveals an interesting parallel between the evolution of physics' recognition of the role symmetry plays in our understanding of the natural world and a similar evolution in mathematics. This parallel structure is stressed not merely to promote an analogy that coaxes us into defining mathematics in terms of symmetry merely because physics does. Rather, we aim to show that philosophical problems about mathematics are resolved when considered in light of development of the role of symmetry. 

Our paper proceeds as follows: in Section 2 we briefly describe the prominent role that symmetry plays in the development of the laws of physics. Section 3 is the core of the paper where we use the concept of symmetry to put forth a novel definition of mathematics and the objects of its study. This definition is employed in Section 4 where we tease out some implications for current questions in the philosophy of mathematics. The paper ends with an appendix that offers a category theoretic formalization of the idea behind reconstructing mathematical structure by looking at the symmetries of a system. 

\section{Physics via Symmetries}

A powerful history of physics can be given in terms of the role symmetry plays in it. There is an ever-expanding place for symmetry in understanding the physical world.\footnote{For philosophical introductions to symmetry see the Introduction in \cite{Brading2003}, \cite{Brading2007}, \cite{Brading2008} and \cite{Bangu2012}. For a popular introduction to the physical issues see \cite{Lederman2004}.} We briefly outline this expanding history from its origins in the classification of crystals to its role in explaining and defining the laws of nature. 

``Symmetry'' was initially employed as it is in everyday language. Bilateral symmetry, for example, is the property an object has if it would look the same when the left and the right sides are swapped. In general an object has symmetry if it appears the same when viewed from different perspectives. A cube thus has six-sided symmetry while a sphere is perfectly symmetrical because it looks the same from any of its infinite sides. 

Pierre Curie formulated one of the earliest symmetry rules about nature. He showed that if a cause has a certain symmetry, then the effect will have a corresponding symmetry. Usually, this rule is stated in its contrapositive form: if an effect has a certain asymmetry, then the cause must also have a certain asymmetry (see e.g. \cite{Stewart92}: 8). Although Curie was mostly concerned with crystals and the forces that create them, the ``Curie Symmetry Principle'' has been employed throughout physics. 

Physicists have generalized the term ``symmetry'' from descriptions of objects to descriptions of laws of nature. A law of nature exhibits a symmetry when it can be viewed from multiple perspectives and still remain the same. We say that such a law is ``invariant'' with respect to some change of perspective. 

An early noticed example of a symmetry exhibited by a law of nature is the fact that the results of an experiment are not altered when the location of an experiment is changed. A ball can be dropped in Pisa or in Princeton and the time needed for the ball to hit the ground will be the same (all other relevant factors being equal). We thus say that the law of gravity is invariant with respect to location. This fact about locations of experiments was so obvious that scientists did not notice or articulate it as a type of symmetry for some time. 

Similarly, the laws of nature are invariant with respect to time. A physical process can be studied today or tomorrow and the results will be the same. The directional orientation of an experiment is also irrelevant and we will get identical results (again, {\em ceteris paribus}) regardless of which way the experiment is facing. 

Formally, we say the equations of classical motion are invariant under translations of spatial coordinates, $x \lra x+\Delta x$, time coordinates, $t \lra t+\Delta t$, and of rotation $r \lra r+ \Delta r$. For each of these transformations, the laws of physics can be transformed to get the same results. 

Classical laws of motion formulated by Galileo and Newton display more complicated symmetries called ``Galilean invariance'' or ``Galilean relativity.'' These symmetries show that the laws of motion are invariant with respect to the inertial frame of reference. This means that the laws of motion remain unchanged if an object is observed while stationary or moving in a uniform, constant velocity. Galileo elegantly illustrates these invariances describing experiments that can be performed from inside a closed ship (\cite{Galilei1953}). Mathematicians have formulated this invariance by studying different transformations of reference frames. A ``Galilean transformation'' is thus a change from one frame of reference to another that differs by only a velocity: $x \lra x+vx$, where $v$ is a constant (small) velocity. Mathematicians realized that all these Galilean transformations form a group which they called the ``Galilean group.'' In broad terms then, the laws of classical physics are invariant with respect to the Galilean group. 

We explain what it means for laws of physics to be invariant with respect to Galilean transformations as follows. Imagine a passenger in a car going a steady 50 miles per hour (neither accelerating nor decelerating) and not turning. The passenger is throwing a ball up and catching it. To the passenger (and any other observer in the moving car) the ball is going straight up and straight down. However to an observer standing still on the sidewalk, the ball leaves the passenger's hand and is caught by the passenger's hand but it does not go straight up and down. Rather the ball travels along a parabola. The two observers are not observing different laws of physics in action. The same law of physics applies to them both. But when the stationary observer sees the ball leaving the passenger's hand, it does not only have a vertical component. Rather, the hand is also moving forward at 50 miles per hour and hence the ball has a horizontal component. The two observers see the phenomena from different perspectives, but the results of the laws must be the same. Each observer must be able to use the same law of physics to calculate where the ball will land despite the fact that they made different observations. Thus, the law is invariant with respect to the ability to swap the two perspectives and still get the same answer. The law is symmetric.  

Another way of saying this is that observers cannot determine whether they are moving at a constant velocity or standing still. The laws of physics cannot be used to differentiate between them because the laws operate identically from  either perspective. 

One of the most significant changes in the role of symmetry in physics was Einstein's formulation of Special Theory of Relativity (STR). When considering the Maxwell equations that describe electromagnetic waves Einstein realized that regardless of the velocity of the frame of reference, the speed of light will always appear to be traveling at the same rate. Einstein went further with this insight and devised the laws of STR by postulating an invariance: the laws are the same even when the frame of reference is moving close to the speed of light. He found the equations by assuming the symmetry. Einstein's insight was to use symmetry considerations to formulate laws of physics. 

Einstein's revolutionary step is worth dwelling upon. Before him, physicists took symmetry to be a property of the laws of physics. It was only with Einstein and STR that symmetries were used to characterize relevant physical laws. The symmetries became {\em a priori} constraints on physical theory. Symmetry in physics thereby went from being an {\em a posteriori} sufficient condition for being a law of nature to an {\em a priori} necessary condition. After Einstein, physicists observed  phenomena and picked out those that remained invariant when the frame of reference was moving close to the speed of light and subsumed it under a law of nature.  In this sense, the physicist acts as a sieve, capturing some invariant phenomena,  calling it a law of physics, and letting the other phenomena go. 

In group theoretic parlance, physicists described ``Lorentz transformations'' - a change from one frame of reference to another that differs by a constant velocity: $x \lra x+vt$ where $v$ is a constant velocity that could be close to the speed of light. All the Lorentzian transformations form the ``Lorentz group.'' It is obvious that the Galilean group is a subgroup of the Lorentz group. The laws of STR are invariant with respect to the Lorentz group.

The General Theory of Relativity (GTR) advanced this idea further by incorporating changes in velocity. From then on the laws of nature were understood as invariant even when acceleration is taken into account. Again, Einstein postulated that the laws of physics should be the same whether or not there are large bodies near the observer and demanded that the laws be invariant if the observer is accelerating. 

In 1918 symmetry became even more relevant to (the philosophy of) physics with Emmy Noether's celebrated theorem. She proved a theorem that connected symmetry to the conservation laws that permeate physics. The theorem states that ``For every continuous symmetry of the laws of physics, there must exist a related conservation law. Furthermore, for every conservation law, there must exist a related continuous symmetry.'' For example, the fact that the laws of physics are invariant with respect to space corresponds to conservation of linear momentum. The law says that within a closed system the total linear momentum will not change and the law is ``mandated'' by the symmetry of space. Time invariance corresponds to conservation of energy, orientation invariance corresponds to conservation of angular momentum, etc (see e.g. \cite{Feynman67} Ch. 4, \cite{Weinberg92} Ch VI, and \cite{Stenger2006} for discussion). Noether's theorem had a profound effect on the working of physics. Whereas physics formerly first looked for conservation laws, it now looked for different types of symmetries and derived the conservation laws from them. 

Currently, particle physics is one of the more interesting places where symmetries are sought and found. Three initial symmetries were postulated: an invariance, $P$, for parity with respect to space reflection (i.e. swapping right for left), the translation, $T$, from $t \lra -t$ (i.e. swapping going forward in time for going backward), and the charge replacement, $C$, of a particle with a corresponding anti-particle. After these, the remainder of the history of particle physics continues the effort to find more symmetries. 

Gauge symmetry is a more abstract symmetry that has played an important role in modern particle physics. Gauge symmetry is what allows a theory that describes different unobservable fields to explain identical observable quantities that emerge from the fields. In short, it is the idea that the laws of physics remain the same no matter how they are described. 

Victor Stenger, a contemporary expositor of symmetry's role in physics, unites the many different types of symmetries under what he calls ``point of view invariance'' (POVI). That is, all the laws of physics must be symmetric with respect to POVI. The laws must remain the same regardless of how they are viewed. Stenger  (\cite{Stenger2006}) demonstrates how much of modern physics can be 
recast as laws that satisfy some type of POVI. 

Symmetry also plays a role in more speculative areas of physics. 
Our best models for going beyond the standard model are attempts to unify all interactions in nature. One of them, supersymmetry, postulates that there is a symmetry that relates matter to forces in nature. Supersymmetry requires us to postulate the existence of a partner matter particle for every known particle that carries a force, and a force particle for every matter particle. The idea here is that the laws of physics are invariant if we swap all the matter for all the force. None of the partner particles have yet been discovered, but because they are mandated by the symmetries that is what scientists are looking for. 

Symmetry, as we have described it, is only part of the story. In numerous cases the laws of physics actually violate a symmetry law and break into different laws via a mechanism known as ``symmetry breaking'' and these broken symmetries are as conceptually important as the symmetries themselves.  The way a symmetry breaks determines certain constants of nature. The question of why a symmetry should break in one way and not another is essentially unanswerable. Researchers are at a loss when they leave the constraints set by symmetry.

Recent excitement over the discovery of the Higgs Boson reveals a triumph of the role of symmetry in physics. Scientists postulated that there was a symmetry in place at the time of the big bang, and it was only when this symmetry was ``broken'' via the ``Higgs mechanism'' that it was possible for mass to exist. By discovering the Higgs Boson physics is able to provide the mechanism by which mass was produced out of the perfect symmetry of the initial state of the universe. The Higgs mechanism was postulated only on the strength of the presumed symmetries.\footnote{See \cite{Bangu2008} for a related discussion of the discovery of the $\Omega-$ particle.}  The recent discovery of the Higgs Boson as the culmination of an extensive research program has further vindicated the methodology of postulating symmetries to discover fundamental properties of the universe. 

Physics also respects another symmetry, which as far as we know has not been articulated as such. The symmetry we refer to is similar to the symmetry of time and place that was obvious for millennia but not articulated until the last century. Namely, a law of physics is applicable to a class of physical objects such that one can exchange one physical object of the appropriate type for another of that type with the law remaining the same. Consider classical mechanics. The laws for classical mechanics work for all medium sized objects not moving close to the speed of light. In other words, if a law works for an apple, the law will also work for the Moon. Quantities like size and distance must be accounted for to get a true law, but when a law is stated in its correct form, all the different possibilities for the physical entities are clear, and the law works for all of them. We shall call this invariance of the laws of nature, {\it symmetry of applicability}, i.e. a law is invariant with respect to the object it is applied to. We shall see in the next section that this is very similar to a type of symmetry that is central to mathematics. 

To sum up our main point, philosophically the change in the role of symmetry has been revolutionary. Physicists have realized that symmetry is the {\em defining property} of laws of physics. In the past, the ``motto'' was that 

\centerline{\underline{The laws of physics respect symmetries.}}

In contrast, the view since Einstein is: 

\centerline{\underline{That which respects symmetries are laws of physics.}}

\noindent In other words, when looking at the physical phenomena, the physicists picks out those phenomena that satisfy certain symmetries and declares those classes of phenomena to be a law of physics. 
Stenger summarizes this view as follows ``\ldots the laws of physics are simply restrictions on the ways physicists may draw the models they use to represent the behavior of matter'' (\cite{Stenger2006}: 8). They are restricted because they must respect symmetries. 
From this perspective, a human physicist observing phenomena is not passively taking in the laws of physics. Rather the observer plays an active role. She looks at all phenomena and picks out those that satisfy the requisite symmetries. 
That is, she takes any extant structure she finds and constructs laws of nature from them.

This account explains the seeming objectivity of the laws of physics. In order for a set of phenomena to be considered a law of physics, it must hold in different places, at different times, be the same from different perspectives, etc. If it does not have this ``universality,'' then it cannot be thought of as a law of physics. Since, by definition, laws of nature must have these invariances,  they appear independent of human perspective. Symmetry thus became fundamental to the philosophical question of the ontology of laws of physics.\footnote{Whether or not there are laws of nature at all or whether they should be eliminated in favor of symmetries in a matter of considerable controversy among philosophers of science. See van Fraassen (\cite{Fraassen89}) and Earman \cite{Earman2004} for stronger and weaker versions of eliminationist views on this issue. Our account is agnostic about this controversy.} 

\section{Mathematics via Symmetries}

Consider the following three examples: 


(1) Many millennia ago someone noticed that if five oranges are combined with seven oranges there will then be twelve oranges in total. Not long thereafter, it was noticed that when five apples are combined with seven apples there is a total of twelve apples. That is, if we substitute apples for oranges the rule remains true. Some time later, in a leap of abstraction, a primitive mathematician formulated a rule that in effect says $5+7 =12$. This last short, abstract statement holds for any objects whatsoever that can be exchanged for oranges. The symbols represent any abstract or real entities such as oranges, apples or manifolds. 


(2) Ancient Egyptians studied different shapes in order to measure ({\it metry}) the earth ({\it Ge}) so that they can have their taxes and inheritance properly assessed. Their drawings on papyrus represented shapes that could be used to divide up the fields on the banks of the Nile. Archimedes could make the same shapes with sticks in the sand. Drawings can accurately describe properties of these shapes regardless of what they represent: plots of land, Mondrian paintings, shipping containers, or whatever. 

In modern times, mathematicians talk about numerous geometrical or topological theorems such as the Jordan Curve Theorem. This statement says that any non-self-intersecting (simple) closed continuous curve in the plane splits the plane into an ``inside'' and an ``outside.'' So for every closed continuous curve there are two regions. If you exchange one curve for another, you will change the two regions. The curve could represent a children's maze or a complicated biological drawing, and the statement remains the same. 



(3) One of the central theorems in algebra is Hilbert's Nullstellensatz. This says that there is a relationship between ideals in polynomial rings and algebraic sets. The point is that for every ideal, there is a related algebraic set and vice versa. In symbols: $$I(V(J))=\sqrt{J}$$ 
for every ideal $J$. If you swap one ideal for another ideal, you get a different algebraic set. If you change the algebraic set, you get a different ideal. This theorem relates the domains of algebra and geometry and is the foundation of algebraic geometry. 

In these three examples we made use of ways of changing the semantics (referent) of mathematical statements. We swapped oranges for apples, changed shapes, transformed curves, and switched ideals. Our central claim is that this ability to alter what a mathematical statement denotes is the fundamental defining property of mathematics. Of course not all transforms are permitted. If we swap some of the oranges for some of the apples, for example, we will not necessarily get the same true mathematical statement. If we substitute a simple closed curve for a non-simple closed curve, the Jordan Curve Theorem will not hold true. Such transformations are not legal. We can only change what the statement means in a structured way. Call this structured changing that is permitted a {\em uniform transformation}. 
Our main point is that this uniform transformation and the fact that statements remain true under such a transformation comprises a type of symmetry. If a statement is true, and we uniformly transform the referent of the statement, the statement remains true. Mathematical statements are invariant with respect to uniform transformations. We call the property of a statement that allows it to be invariant under a change of referent {\em symmetry of semantics}. The statement remains the same despite the change of semantic content. 

Every mathematical statement defines a class of entities which we call its ``domain of discourse.'' This domain contains the entities for which the uniform transformation can occur. When a mathematician says ``For any integer $n \ldots$,'' ``Take a Hausdorff space \ldots'', or ``Let $C$ be a cocommutative coassociative coalgebra with a involution \ldots'' she is defining a domain of discourse. Furthermore, any statement that is true for some element in that domain of discourse is true for any other. Notice that the domain of discourse need not be only a single class of entity. Each statement might have $n$-tuples of entities, like an algebraically closed field, a polynomial ring and an ideal of that ring. Every mathematical statement has an associated domain of discourse which defines the entities that we can uniformly transform.


Different domains of discourse are indicative of different branches of mathematics. Logic deals with the classes of propositions while topology deals with various subclasses of topological spaces. The theorems of algebraic topology deal with domains of discourses within topological spaces {\em and} algebraic structures. One can (perhaps naively) say that the difference between applied mathematics and pure mathematics is that, in general, the domains of discourse for applied mathematical statements are usually concrete kinds of entities while the domains of discourse for pure mathematical statements are generally abstract entities. 

From this perspective, one can see how variables are so central to mathematical discourse and why mathematicians from Felix Klein to Tarski to Whitehead touted their import for mathematics (\cite{Epp2011}). Variables are placeholders that tell how to uniformly transform referents in statements. Essentially, a variable indicates the type of object that is being operated on within the theory and the way to change its value within the statement. For example in the statement $$a \times (b + c) = (a \times b)+(a\times c)$$ which expresses the fact that multiplication distributes over addition, the $a$ shows up twice on the right side of the equivalence. If we substitute something for $a$ on the left, then, in order to keep the statement true, that substitution will have to be made twice on the right side. In contrast to $a$, the $b$ and $c$ each occur once on each side of the equation. Transformations to those entities will be simpler. Again, the variables show us how to uniformly transform the entities. 

For us, a mathematical object is any entity that can be used in a mathematical statement. So oranges, apples, and stick drawings in the sand are mathematical objects. As long as we can transform those objects into other objects within the same domain of discourse, then they are mathematical objects. We can transform seven oranges into the elements of the set ${7} = \{0,1,2,3,4,5,6 \}$ and we give equal status to each of them as mathematical objects. Mathematicians prefer to use $7$ because of the generality it connotes. But representing it this way  is misleading because seven oranges are just as good at representing that number in the sentence $5+7=12$. Any statement about the number seven can be made with a transformation of the elements from the set of seven oranges as well as the set $7$.\footnote{Frege's influence on this definition should be evident.  Consider the collection of finite sets. A finite number for Frege consists of the equivalence class of the finite sets where two sets are equivalent if there is an isomorphism from one set to another. When we talk of the equivalence class 5 we are ignorant of the element of that equivalence class. Are we talking about 5 apples or 5 cars? Here we are talking about the equivalence class as opposed to the group (groupoid) but our account generalizes this. We equally account for not only the finite numbers but all objects and all operations of mathematics.} (As we shall explore later, this accounts for practices in applied mathematics.) The mathematical statements made by applied mathematicians are no less true than the statements made by pure mathematicians. 

At this point  two  objections can be raised regarding our stance on metaphysics. First, given what we usually mean by ``apple'' and what we usually mean by ``mathematical object'' an apple is not a mathematical object, so how can we substitute apples for $7$; isn't this a kind of category error? Our answer is that it is not, because by our definition a mathematical object is merely an object that is amenable to mathematical treatment. We can assert this because it is the only distinction that matters to mathematics and mathematicians. This affords us the philosophical advantage of preserving a uniform semantics across all kinds of objects which our ontology countenances, and our ontology stays neutral about what it means for $7$ to exist. For $7$ to exist means nothing above the fact that it can be swapped in a uniform transformation with other objects (or ``objects'') in a way that preserves symmetry of semantics. 

Secondly, Doesn't allowing anything to be a mathematical object obscure the question of what mathematical objects are? This, we say, is beside the point. Once we account for why and how mathematics works in all contexts in which it does work, demanding an account of what mathematics {\it really} deals with, is outside the scope of what mathematics could and should aim for. We have achieved a mathematical analogue of empirical adequacy.  Looking for the underlying homunculus or soul or even constituent parts of the numbers, manifolds, sets, infinitesimals, or whatever advances neither mathematics nor its philosophy. 

Symmetry of semantics should look familiar in this context to logicians and model theorists as the definition of validity. A logical formula is valid if it is true under every interpretation. That is, it must be true for any object in the domain of discourse. Thus, symmetry of semantics is not a radical idea. Rather, the novelty is viewing validity as a type of symmetry. We shall see that this symmetry is as fundamental to mathematics as many symmetries are to physics. 

While it should now be obvious that mathematics satisfies symmetry of semantics, we would like to do for mathematics what is now natural in physics. Rather than understanding mathematics as a discipline that satisfies symmetry of semantics, we define mathematics as that which satisfies these symmetries. In other words, whereas we used to understand that: 

\centerline{\underline {A mathematical statement satisfies symmetry of semantics}.}

\noindent we replace this by making symmetry of semantics the defining property of mathematics; so:

\centerline{\underline {A statement that respects symmetry of semantics is a mathematical statement}.}

\noindent In other words, given the many expressible statements a mathematician finds, her job is to distill and organize those that satisfy symmetry of semantics. In contrast, if a statement is true in one instance but false in another instance, then it is not mathematics. In the same way that the physicist acts as a ``sieve'' and chooses those phenomena that satisfy the required symmetries to codify into physical law, so too the mathematician chooses those statements that satisfy symmetry of semantics and dubs it mathematics. 

Moreover, just as the physicist constructs laws of nature to subsume wide varieties of phenomena under the rubric of one law, so to does the mathematician collect many instantiations under one mathematical statement.  

Many statements in general do not satisfy symmetry of semantics. Poetry and song are full of lines where one cannot substitute any one word for another. Even some mathematical-type statements like ``If x is like y and y is like z, then x is like z'' simply fail in most cases. The reason is that the word ``like'' is vague and not exact enough to be part of mathematics. 

One may want to object to this view by saying that it is too inclusive: many branches of the sciences and statements in ordinary language satisfy symmetry of semantics, thereby making them part of mathematics. We do not disagree. In fact, many branches of science can be seen as applied mathematics. One would be hard pressed to find a good dividing line between theoretical physics and mathematics, computer science and mathematics, etc.\footnote{Whereas we consider these fields parts of mathematics, Steiner (\cite{Steiner2005}) treats them as applications of mathematics.} Moreover, many ordinary language statements such as ``all women are mortal'' satisfy a certain symmetry of semantics and hence is also part of mathematics.  In both cases this is because mathematics can model many parts of the real world. Many statements in both pure and applied sciences do satisfy symmetry of semantics.  The inability to clearly draw boundaries around these fields is a direct result of the fact that all mathematics respects the same kind of symmetry.\footnote{See also \cite{Longo2013} for an application to theoretical biology.} 

Symmetry of semantics is not the only symmetry that mathematical statements satisfy (though some symmetries seem so obvious that even mentioning them seems strange) and mathematics also satisfies many symmetries that physical phenomena satisfy. For example, mathematical truths are invariant with respect to time: if they are true now, then they will also be true tomorrow. They are invariant with respect to space: a mathematical statement that is true in Manhattan, will also be true on Mars. Other symmetries that mathematics satisfies also usually go unsaid. For example, it is irrelevant who proves a theorem or in what language the proof is done. 

Another symmetry that mathematics has we call {\em symmetry of syntax}. This says that any mathematical object can be described (syntax) in many different ways. For example we can write 6 as $2 \times 3$ or $2+2+2$ or $54/9$. The number $\pi$ can be expressed as $\pi=C/d$, $\pi = 2i\log{1 - i \over 1 +i}$, or the continued fraction $$\pi=3+\textstyle \frac{1}{7+\textstyle \frac{1}{15+\textstyle \frac{1}{1+\textstyle \frac{1}{292+\textstyle \frac{1}{1+\textstyle \frac{1}{1+\textstyle \frac{1}{1+\ddots}}}}}}}$$

Similarly we talk about a ``non-self-intersecting continuous loop,'' ``a simple closed curve,'' or ``a Jordan curve'' and mean the same thing. The point is that the results of the mathematics will be the same regardless of the syntax we use. Mathematicians often aim to use the simplest syntax possible, so they may write ``6'' or $\pi$ instead of some equivalent statement, but ultimately the choice is one of convenience, as long as each option is expressing the same thing. 

Similar to what we saw in the case of physics, by selecting only those statements that are invariant with regard to what the statement is referring to, the mathematician is ensuring that the statement is objective and universal. Since mathematical statements can, by definition, refer to so many different entities that seem to be independent of human beings, they appear to have an independent existence. 

Our definition of mathematics also demystifies why mathematics is so useful to the natural sciences. As we saw in the last section, the laws of physics are invariant with respect to the symmetry of applicability. This means that the laws can apply to many different physical entities. Symmetry of applicability is a type of symmetry of semantics. In detail, symmetry of applicability says that a law of nature can apply to many different physical entities of the same type. Symmetry of semantics says that a mathematical statement can refer to many different entities in the same domain of discourse. When a physicist is trying to formulate a law of physics, she will, no doubt, use the language of mathematics to express this law, because she wants the law to be as broad as possible. Mathematics shares this broadness. The fact that some of the mathematics could have been formulated long before the law of physics is discovered is not so strange. Both the mathematician and the physicists chose their statements to be applicable in many different contexts. The ``mystery'' of the unreasonable usefulness of mathematics melts away.

\section{Some Philosophical Consequences}

So far we have used symmetry and invariance to provide a definition of mathematics. But what questions about the nature of mathematics are thereby answered? 

As we mentioned in the Introduction, our characterization faces the same challenges as any other philosophical account of mathematics: it should provide a coherent ontology and methodology, and a semantics that meshes with mathematical practice. It should also provide a semantics that treats mathematics similar to comparable modes of discourse (or explain why it does not) and provide an epistemology that accounts for why we find mathematical arguments so compelling. Our responses to these challenges diverge significantly from other accounts of the nature of mathematics. 

Our account takes the idea that mathematics determines what objects and statements are mathematical, modifies it, fleshes it out, and embraces the notion that it is those statements which respect symmetry of semantics that determine what  mathematics is. Doing so implies the following: (1) Contrary to an important Platonist assumption (\cite{sep-platonism-mathematics}), mathematics is not independent of humans. Rather, mathematics respects the symmetries humans deem important. (2) Contrary to standard nominalist assumptions, mathematical objects do exist. Anything, including concrete objects, which can be subsumed within a uniform transformation of the proper kind is a mathematical object. (3) Contrary to most accounts of mathematics (especially Platonism), mathematical objects are not ontologically uniform; mathematical objects can be concrete, abstract, or perhaps something else (e.g. theological, fictional, etc?). 

In addition, our account treats mathematical objects semantically no differently from scientific objects. By our definition, a mathematical object is mathematical if it respects the same kinds of symmetries as physical objects respect in physics. So when we say that 8 is larger than 3 and when we say an atom is larger than a proton we are not operating with different semantics for the two statements. They are both true, and they are both true in the same way. Assigning identical truth conditions is sometimes seen as a problem for some kinds of nominalists who have to account for the fact that the latter statement is true, while the former statement is not (because ``8'' and ``3'' do not refer). For us however, the very fact that the sentences assume the same form shows that they have symmetry of semantics. They can be uniformly transformed in the context of a ``$>$'' and {\em that} is what makes them both mathematical statements with the same semantic strategy. 

Given (1) we can also account for why it is generally easy to agree on which bits of mathematics are true or why mathematics looks ``objective'' without having to appeal either to an innate Platonic sense organ  (``Radio Plato'', as Jody Azzouni memorably put it (\cite{Azzouni2008})) that detects mathematics on the one hand, or deny the truth of mathematics on the other. On our account there is no mathematical structure ``out there'' to Platonically ``see.'' The structures come from the symmetries we presuppose.  Given (2) our account of mathematics is also not an ``error theory.'' An error theory is one that commits the speaker to a mistake every time she utters a mathematical statement because mathematical discourse does not actually refer to the objects we seem to commit to when uttering mathematical statements. To the contrary, our account finds it unsurprising that we all agree on mathematical truths because mathematics starts with an {\em a priori} conception that mathematics will obey symmetry of semantics. Given those symmetries, mathematics has no other way to turn out. 

Thus symmetry of semantics makes mathematics look the same to everyone. More precisely, since it is that which has symmetry of semantics that we call ``mathematics,'' our account also accounts for apparent mathematical objectivity, which not all other conceptions of mathematics can do.\footnote{Nozick's \cite{Nozick01} is mostly dedicated to the more general case for the equivalence of invariance and objectivity though he does not address our point.} 

Implications (2) and (3) - that mathematical objects exist and that they are ontologically heterogeneous - also provide our account with the advantage of a broad ontological neutrality. Meaning, since we have no preconceived notion of what numbers are, we need no account of some other equally abstract object and an account of its relation to familiar numbers. Numbers can be whatever we need them to be.  Yes, mathematics can pull off such reductions, going from geometry to algebra or from numbers to logic and set theory, but as we shall see below, that is only a result of the underlying symmetries in mathematics. Our account instead preserves standard mathematical discourse as it can talk about any ``gross'' mathematical object without worrying about what it {\em really} is. That is, we do not say that when mathematicians talk about numbers they are really talking about structures or sets. We say instead that when mathematicians talk about numbers, they are talking about numbers and when they talk about apples, they are talking about apples, even if they talk about a number of them.  This is comparable to the way physics can talk about planets without worrying about their constituent particles. Planets may reduce to strings (if your ontology goes that way) but planetary dynamics is no less physics than M-theory, and it would be odd to suggest that both the laws and the ontology of planets are really built up from the laws and ontology of strings, even if you believe in laws of nature, reduction and bridge laws. 


{\bf Wigner's problem and naturalism.} A philosophical naturalist's interest in the philosophy of mathematics is the alignment of the ontology, epistemology, and especially methodology of mathematics with those established for science. To the extent that we are convinced that the underlying methodology in the natural sciences (or at least in fundamental physics) has been a shift in our understanding of the laws of nature to accommodate an {\em a prioristic} concept of symmetry, then we have certainly put forth a naturalistic treatment of mathematics. The shift in the way we understand mathematics is now aligned with the shift in the way we understand the role of symmetry in physics. 

This treatment explains the problem of the unreasonable effectiveness of mathematics in the natural sciences. This problem as famously articulated by Wigner (\cite{Wigner60}) (on one interpretation) finds the physical science that we discover shockingly related to the mathematics we need to understand it. Yet, every time science needs to articulate a physical concept, it finds that it can turn to mathematics for help; the mathematics will be there. Mark Steiner (\cite{Steiner95}: 154) sees one version of the problem as stemming from the apparent mismatch of methodologies. How can problems emerging from physics be articulated, and even solved, using methods that were designed for a completely unrelated purpose? 

A. Zee, completely independent of our concerns, has re-described the problem as the question of ``the unreasonable effectiveness of symmetry considerations in understanding Nature.'' Though our notions of symmetry differ, he comes closest to articulating the way we approach Wigner's problem when he writes that ``Symmetry and mathematics are closely intertwined. Structures heavy with symmetries would also naturally be rich in mathematics'' (\cite{Zee1990}: 319). 

To reiterate what we said at the end of the previous section, on our account, physics discovers some phenomenon and seeks to create a law of nature that subsumes the behavior of that phenomenon. The law must not only encompass the phenomenon but also a wide range of phenomena. The range of phenomena that is encompassed expresses a symmetry of applicability. So a law must be deliberately designed with symmetry of applicability. Mathematics has a built in ability to express these symmetries because the symmetry of applicability in physics is actually just a subset of the symmetries that mathematics has - i.e. the symmetry of semantics. There is then nothing surprising about the fact that there is some mathematics that is applicable to physics, as the symmetries of physics are a subset of the symmetries of mathematics. Wigner's problem thereby loses its force. 




{\bf Foundations} Our way of understanding mathematics also describes the way we take mathematics to have foundations. 
Typically ontological, methodological, and epistemological features of philosophical programs that address the nature of mathematics describe the way in which mathematics is ``grounded'' in a foundation, if indeed it is so grounded. Ontological foundations are supposed to describe what the constituent nature of mathematical objects ``really are;'' methodological foundations describe the (and the only) methodology that is acceptable throughout mathematics; while epistemic foundations are tasked with explaining why mathematical statements are so convincing as compared to other areas of knowledge (see \cite{Azzouni2005}). 

By now our epistemic story should be clear: Our confidence in the certainty of mathematics has its origins in our {\em a prioristic} concept of symmetry instantiated in mathematics. The fact is, we are certain about mathematical results because we have decided that the only mathematically tractable entities are those entities that are amenable to uniform transformations. 

Our account also explains the intuitions behind reductionist accounts of mathematics, such as those that build numbers out of sets. When we speak of such ontological reductions we are treating ``foundation'' in connection with mathematics as the ability to show that large parts (or even all) of mathematics can be phrased in some system, and that system is ``primitive''.  Commonly, since many parts of mathematics can be reduced to set theory and logic and sets provide such a convenient domain of discourse in everyday settings, set theory is taken to be a good candidate for such a foundational system. 

This is presumably analogous to the conception of fundamental physics that seeks out the particles in which we, in theory, can express our fundamental ontological statements. This search for fundamental laws or fundamental particles is an important part of contemporary physics. But as we have shown, physics has largely abandoned the idea that these programs are ``starting points'' for scientific research. Starting points have been replaced by the presumption that invariances are (methodologically) fundamental. This assumption in turn allows fundamental particles to be found. Programs regarding the foundations of mathematics began by confusing reduction and invariance in a very particular way. 

This can be seen as a two part problem. First, given our ontological conclusions about mathematics - that mathematics is more ontologically similar to science than is usually supposed -  the usual rules of science ought to apply. 
When we look for the ``starting points of mathematics'' it makes little sense for us to look for its fundamental pieces, as we have concluded that mathematics admits a variety of types of objects that have no reasonable expectation of reducing one to another (and in many cases it would anyway be unclear what reduces to what). Mathematics, in other words, has no fundamental parts, nor can we necessarily find one kind of object with which to phrase all the others. The vocabulary of uniform transformations (i.e. the vocabulary of mathematical methodology) is the only way to talk about both abstract and concrete objects. 

Secondly, although typical mathematics can be reduced to sets, sets do not exhibit the correct kind of expressive power or display the right kind of symmetries in mathematics to be the fundamental ``ground''.  Sets only display the symmetry of mathematical objects. That is, set theory shows that all mathematical objects are the same in one way: they can all be ``reduced'' to the same thing. Since all mathematical objects are the same, there are ways in which they can be treated similarly (akin to showing that all (non-fundamental) physical objects reduce to fundamental particles, it fails to make sense of all the other symmetries in nature). As the Ernie and Johnny story (\cite{Benacerraf65}) shows, set theory, for example, does not actually make sense of the myriad of ``fundamental mathematical properties''. The story shows that there are an infinite number of set-theoretical reductions because the number 3, for example, can be reduced to (identified with) $\{\{\{\varnothing\}\}\}$ or $\{\varnothing, \{\varnothing\}, \{\varnothing, \{\varnothing\}\}\}$ or any other way of configuring set-theoretic ``three-ness.'' And this very concern is a symptom of confusing reduction and invariance, not a problem with our understanding of the real nature of numbers. The moral Benacerraf draws is that a number cannot be a set because there is no one set that corresponds to each number. However,  Benacerraf's counterexamples are exactly our point. The fact that we can swap $\{\varnothing, \{\varnothing\}, \{\varnothing, \{\varnothing\}\}\}$ for $\{\{\{\varnothing\}\}\}$ shows that set theory exhibits symmetry of semantics. But it says nothing about how sets or set theory are foundational. The fact that we can swap one set for another and understand both as 3 is not only unsurprising on our account, but expected, because set theory is just another branch of mathematics that has symmetry of semantics, like all the others. 

In part, the real reason set theory does not show anything fundamental about numbers is that it does not account for the way we actually  take mathematics to exhibit invariances. Namely, we take mathematical objects to be invariant in a way that the objects stay the same under a wide range of rule transformations; not just object transformations. 
Therefore, the reason that set theory initially appears intuitively like a ground for mathematics, but nonetheless fails, is because set theory can do one thing that we expect of a physical reduction, namely exhibit an ontological-type reduction of some mathematical objects, but the reduction is methodologically inadequate as it cannot capture what is really important (what is really mathematical) about mathematics. 

As in physics, the fundamental nature of the objects is interesting, but the methodological heuristic of symmetry is what is common to all mathematics. The question should not be ``what are numbers?'' but rather ``what is mathematics?''  Foundationally, methodology, not ontology is important, and the methodology we actually employ is symmetry.\footnote{Our way of understanding symmetry in mathematics is similar to the way Otavio Bueno (\cite{Bueno2006}) understands symmetry in quantum mechanics - as a methodological heuristic.} Normally, however, when we think about foundations in methodological terms we think of the Euclidean model where mathematics requires an axiomatic system plus a system of deduction that allows for the generation or construction of the rest of mathematics (or its contemporary model theoretic version.) But as Azzouni (\cite{Azzouni2005}) has argued this is untenable both because of the G\"{o}del problem of establishing proper axioms and because terms like ``constitute'' or ``construct'' are used metaphorically in mathematics and not literally, as there is no clear ontological hierarchy in mathematics. Symmetry is a deeper methodological foundation neither countenancing an axiomatic nor ontological hierarchy, and thereby skirting Azzouni's worries.  Thus,  we offer a new conception of mathematics which  we claim is foundational in the sense that it describes the mathematical methodology, and accounts for the way mathematics behaves and the way mathematicians treat mathematics in practice.  

{\bf Mathematical practice}: Our considerations about symmetry as the underlying nature of mathematics emerge directly from consideration of both the large and small scale aspects of mathematical practice. We see it as a {\it sine qua non} of any conception of mathematics that it square with the way contemporary mathematics works. If mathematicians cannot (be made to) see their craft in an approach to mathematics, so much the worse for the approach. 

The day-to-day job of the mathematician is proving theorems. Contrary to the impression given by the typical finished mathematics paper, mathematicians do not generally posit a theorem and then proceed to prove it from axioms. In reality a mathematician has an intuition and formulates some statement. The mathematician tries to prove this statement but almost inevitably finds a counterexample. A counterexample is a breaking (violation) of the symmetry of semantics; there is some element in the supposed domain of discourse for which the statement fails to be true. The mathematician then proceeds to restrict the domain of discourse so that such counterexamples are avoided. Again our indefatigable mathematician tries to prove the theorem but fails, so she weakens the statement. Iterating these procedures over and over eventually leads to a proven theorem. The final theorem may only vaguely resemble the original theorem the mathematician wanted to prove. In some sense rather than saying that the ``proof comes to the theorem,'' we might say that ``the theorem meets the proof half way.'' The mathematician acting as a ``sieve'' sorts out those statements that satisfy symmetry of semantics from those that do not, and only those that satisfy this symmetry are reported in the circulated and published paper.\footnote{It should come as no surprise that mathematics (or science) does not emerge from the mathematician as cleanly as it appears to from the mathematical paper. The classic (though itself quite idealized) discussion of actual mechanics of mathematical proof can be found in Polya's {\em How to Solve It} (\cite{Polya57}).} Lakatos' ``rational reconstruction'' in his \cite{Lakatos76} may be cited as an example of this constant struggle to preserve the symmetry of semantics of the Euler formula. 

As with physics, we observe symmetry on various levels. It is not just with individual proofs, but also on the level of mathematical programs that we see the influence of the notion of symmetry. 

Discovering symmetries in broad fields of mathematics is like discovering symmetries in fields of physics. Discovering symmetries in physics relies on the idea that there are domains in which transformations are allowed. Physics progresses by ever larger numbers of allowable phenomena getting subsumed under a given domain as when Maxwell unified electrical theory and magnetism or when terrestrial and plantary mechanics were united by Newton. Similarly in mathematics symmetries are discovered when we find that seemingly different mathematical phenomena are really in the same category as an already known transformation and are thereby subsumed under a larger domain; we discover that a new larger class of entities can be uniformly transformed. In other words, we find that there is a union of different domains of discourses which were previously assumed to be comprised of non-interchangeable entities. The ``monstrous moonshine'' conjecture is one famous case of such unification. In the late 1970's John McKay noticed a completely unexpected relationship between the seemingly different areas of the ``monster group'' and modular functions. Legend has it that when McKay first heard that the number 196,884 appears in both areas, he shouted ``moonshine'' as a term of disbelief. Deeper connections between the monster group and modular functions have since been shown. 

Another example of such unification is algebraic topology. Researchers realized that there is a certain similarity when one takes maps between two topological spaces into account and when one takes homomorphism between two groups into account. That is, there is a relationship between topological phenomena and algebraic phenomena. Mathematicians went on to use this similarity to try to classify certain topological structures. Category theory grew out of this unification and essentially became a tool for much more unification. By being ``general abstract nonsense,'' as category theory has been called, that is about ``nothing,'' it can thereby be about everything; and hence its language can be used in many different areas of mathematics. As with symmetry, such unification advances mathematics by giving mathematicians an opportunity to discover more general theorems with wider applications and allows them to apply techniques from one domain to the other.\footnote{Philip Kitcher (e.g. \cite{Kitcher76} and \cite{Kitcher89}) touts the importance of these types of cases for mathematics and uses them in the service of demonstrating the existence of mathematical explanation. He argues that the act of unification in a very specific way is how scientists and mathematicians explain. There is some important sense in which we agree with what Kitcher does, though we draw very different lessons than he does.  A discussion of Kitcher, unification, and category theory can be found in \cite{Zelcer09}. Emily Grosholz has studied domain unifications extensively. See e.g. \cite{Grosholz00}, \cite{Grosholz80}, and \cite{Grosholz85}.} 

The final way we see the impact of symmetry in mathematical practice is by analogy with symmetry breaking. An analysis of symmetry breaking in physics can provide a description of the nature of physical phenomena. Physical constants, for example, are what happened when symmetries broke the way they did. Symmetry breaking is also an important part of mathematical practice. Probably the first example of a mathematical broken symmetry was discovered by Hippasus, a student of Pythagoras. The Pythagoreans believed that every number is rational. Hippasus showed that the diagonal of a square has length $\sqrt{2}$ but that it is not a rational number. The idea that {\it every} number is rational was thrown overboard (together with Hippasus). That is, $\sqrt{2}$ was the first element in the domain of discourse known to the Greeks as numbers that showed that the domain must be split or broken in two. This discovery begins the long fruitful history of rational and irrational numbers. For another example consider the many problems in computability theory. This area of theoretical computer science was started by Turing and others in the 1930s. In the 1960s researchers realized that although there are many problems that are decidable/solvable by a computer, there are some problems that take an exponentially long time to solve. The Euler Cycle Problem asks to find a cycle in a given graph that hits every edge exactly once. In contrast, the Hamiltonian Cycle Problem asks to find a cycle in a given graph that hits every {\it vertex} exactly once. Whereas there is a nice polynomial algorithm to solve the Euler Cycle Problem, there is no known polynomial algorithm for the Hamiltonian Cycle Problem. This realization - that the usual methods of solving computational problems fail sometimes, created the entire important field of computational complexity theory. 

Thus there are various ways in which symmetry considerations aptly describe mathematical practice in the same way they describe scientific practice. 

{\bf Beauty.} A final note. Ever since Plato's {\em Timaeus} it has not been uncommon for philosophers, mathematicians, and scientists to praise proportion and symmetry in mathematical laws and equations.\footnote{See Steiner's \cite{Steiner98} for a more thorough discussion of beauty in mathematics. There he argues that an anthropocentric concept of beauty plays an important role in mathematical practice.} Symmetry is often said to be a part of the beauty of mathematics or the natural laws in the same way that Aristotle's {\em Poetics} sees beauty in the symmetry and proportion of the composition of artistic objects. By understanding the nature of mathematics as that which matches our notions of symmetry, as we have already firmly entrenched in our understanding of the natural world, we account for the beauty that many mathematicians see in mathematics. 

It is hard to think of a discussion of beauty - mathematical or otherwise - that either does not  take symmetry into account or apologize for not doing so. But one important lesson we can learn from the above discussion of the symmetrical nature of mathematics and the nature of mathematical symmetry is that thinking of mathematics as beautiful is getting it backward. In actuality, what we are saying is that which is beautiful, is mathematics. 



\section{Appendix: Algebraic Structure via Symmetries in Category Theory}
�\begin{flushright}{\it Thus we seem to have partially \\ demonstrated that even in foundations, \\ not substance but invariant form is the carrier\\ of the relevant mathematical information.}\\
F. William Lawvere \\ 
``An Elementary Theory of the Category of Sets'' \cite{LawvereFW:eletcs}
\end{flushright}

In this paper we have highlighted that you can find mathematical laws by looking at the symmetries that those laws respect. In a sense, we are ``reconstructing'' the mathematical laws by looking at the symmetries of the laws. It turns out that algebraists and category theorists have been doing this for a while in their own contexts. In fact, this was the inspiration and impetus for some of the ideas in this paper. Here we describe these categorical ``reconstructions'' and their meanings.\footnote{For a thorough discussion of the role of symmetry in category theory see \cite{marquis2009geometrical}.}  

In effect, what this appendix will show is that one can think of uniform transformations as something carried out by homomorphism between algebraic structures. Saying that mathematics is what respects uniform transformations is saying that algebraic operations are preserved by homomorphisms. In essence there is a duality between operations and homomorphism as there is a duality between mathematics and uniform transformations. 

As a broad outline, this appendix will explain the following ideas. In 1963 F. William Lawvere defined a type of category called an {\it algebraic theory}. This is a formal categorical way of describing a type of algebraic structure. For example, there will be a category that is the algebraic theory of groups which will be a description of what it means to be a group. There will also be an algebraic theory of rings, an algebraic theory of monoids, an algebraic theory of semigroups, etc. For each algebraic theory, one can construct a category of models\footnote{It should be noted that mathematicians and physicists use the word ``model'' in two exactly opposite ways. For a mathematician a model, as in model theory, is an instantiation of an abstract ideal. There is some theory and a model of that theory is a representation of the ideal theory. A model is some concrete version of an ideal. In contrast, a physicist uses the word as an abstract ideal. There are physical phenomena and they seek to ``model'' them. To physicists a model is an abstract ideal that makes sense of the many concrete phenomena.} for that theory. So there will be a category of groups, a category of rings, a category of monoids, etc. The more interesting aspect to us will be that Lawvere went on show how to go the other way: for a given category of algebraic structures, one can ``reconstruct'' or select its algebraic theory. In other words, one can look at the models of some theory (and its symmetries) and describe its algebraic theory. This is analogous to the way physicists reconstruct or select the laws of physics by looking at the phenomena of the physical world. This is also related to our focus on reconstruction and selection of the mathematical laws by looking at those laws that satisfy symmetry of semantics. 

This appendix assumes some basic category theory. The simple notions of a category, functor, and natural transformation will be needed. For simplicity sake, a category is a ``souped-up'' directed graph with a multiplication of arrows. That is, given two arrows where the source of one is the same as the target of the other, there is an arrow that combines these two arrows. A functor between two categories is simply a ``souped-up'' graph homomorphism. That is, it is a function that takes objects of one category to objects of another category and arrows of one category to arrows of the other. The functor should respect the composition of morphisms. A natural transformation is a type of map from one functor to another functor. With these definitions in mind, one should be able to understand what follows. 

A category is often thought of as some type of collection of different objects that have some structure in common (like the category of topological spaces, the category of sets, the category of rings, etc) we can also think of a category as a structure itself. For example, category theorists think of a partial order as a category. Also, one can think of a single group as a category with only one object and all the morphisms invertible. Lawvere showed how to think of a description of algebraic structure (in logic and model theory this is called a ``signature''; in universal algebra this is called a ``clone'') as a category. An algebraic theory is a category whose objects are the natural numbers and directed edges $f:m \lra n$ correspond to an algebraic operation that accepts $m$ objects and outputs $n$ objects. In particular, a directed edge $f:m \lra 1$ corresponds to the classical notion of an $m$-ary operation. A directed edge $f:0 \lra 1$ corresponds to an operation that has no inputs but just has a particular output. This corresponds to a constant of the structure, that is, a particular element of the structure. Composition of edges correspond to composition of different algebraic operations. Two edges that are equal to each other in an algebraic theory correspond to an identity of the algebraic structure. So for example the identity $a*(b+c)=(a*b)+(a*c)$ will correspond to the fact that two different directed edges from 3 to 1 will be set equal to each other. 

There are many examples of algebraic theories. Lets look at the algebraic theory for groups, $T_{group}$. This category will have the following: 
\begin{itemize} 
\item a morphism $+:2 \lra1$ which correspond to the group multiplication operation, 
\item a morphism $-:1 \lra 1$ which corresponds to the inverse operation, and 
\item a morphism $e:0 \lra 1$ which will correspond to the identity element of the group.
\end{itemize}
These morphisms generate the other morphisms in
$T_{group}$. For example there is a unique
$+$ map from $3$ to $1$ that is the extension of the associative addition: $+(id \times +):3 \lra 1$ 

One can visualize a
small part of this algebraic theory as 
$$
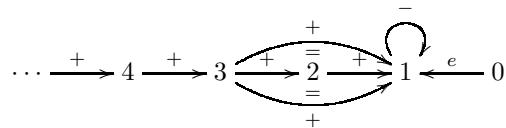
\xymatrix{
\cdots \ar[r]^+ & 4 \ar[r]^+ & 3 \ar[r]^+\ar@/_1pc/[rr]_+^= 
\ar@/^1pc/[rr]_=^+ & 
2 \ar[r]^+ & 1 \ar@(ul,ur)[]^{-} & 0\ar[l]_e
}$$
where the two morphisms from 3 to 1 are the two (equal) ways of
associating three objects.
In some 
sense, $T_{group}$ is the ``shape'' or ``template'' of groups. 

The morphisms must satisfy the following identities:
$$\xymatrix{
3 \ar[dd]_{+ \times id}\ar[rr]^{id \times +}&&2\ar[dd]^+
&& 2\ar[dd]_{-( ) \times id}&&1\ar[rr]^\Delta \ar[ll]_\Delta
\ar[d]_{!} &&2\ar[dd]^{id \times -( )}\\
&& && & &0\ar[d]_e \\
2\ar[rr]_+&&1&& 2\ar[rr]_+&&1&&2\ar[ll]^+
}$$
$$\xymatrix{
1\times 0 \ar[rr]^{id \times e} \ar[rrd]_\simeq && 2 \ar[d]_+ &&
0 \times 1 \ar[ll]_{e \times id} \ar[lld]^\simeq \\
&& 1}$$
where $\Delta$ is the usual ``diagonal'' map and $!$ is the unique map
from $1$ to $0$.

Other examples are the algebraic theory of monoids, $T_{monoid}$. This category will not have the $-:1 \lra 1$. The algebraic theory of rings, $T_{ring}$, will have what a group has but will also have a morphism $*:2 \lra 1$ which will correspond to the second operation that a ring must have. For every ring $R$, there will be an algebraic theory, $T_{R-mod}$, of $R$ modules. They will be the same operations as the theory of groups but with the added operations $r:1 \lra 1$ for each $r \in R$. That is, an $R$-module is simply a group with an action of the ring R on it. For most algebraic structures there will be an algebraic theory that describes this structure. It pays to mention that there does not exist an algebraic theory for fields simply because the division operation $2 \lra 1$ is not a total operation. It is not defined when the second element is 0. Such a theory is called ``essentially algebraic'' and can be described by generalizations of algebraic theories called ``sketches'' (\cite{BarrM:toptt}). 

For a given algebraic theory $T$, an algebra in sets is a functor $A$ from $T$ to $\Set$ that takes $1$ to some set $A(1)$ and takes $n$ to the $n$-th product of the set $A(1)$. That is, an algebra $A$ is a functor $A:T \lra \Set$ such that 
$$A(n) = A(1)\times A(1)\times \cdots \times A(1).$$
This ensures that $A$ will take a morphism $n \lra 1$ to an operation that goes from $n$ copies of some set $A(1)$ to the single copy of $A(1)$. Identities between two operations will correspond to identities of the set operations. Let us examine carefully the concept of a group. An algebra for the theory of a group is a functor $A:T_{group} \lra \Set$ that satisfies the above requirement. The underlying set of the group will be $A(1)$ and the main group operation will be
\begin{itemize} 
\item The group multiplication operation: $A(+):A(2)=A(1)\times A(1) \lra A(1).$ 
\item The group inverse: $A(-):A(1) \lra A(1).$ 
\item The identity element: $A(e):A(0)=\{\star \} \lra A(1).$
\end{itemize}
In a sense, an algebraic theory is like a cookie cutter that is used to make many different cookies. The right shape and proper functors from the algebraic theory to the category of sets picks out a structure of the shape.
One can take algebras in other categories besides $\Set$. All that is needed is a category with products. For a category $\C$, an algebra is a functor $A:T \lra \C$ that respects the product. In particular, the category of topological groups is the category algebras for a the theory of groups, $T_{group}$, in ${\bf Top}$, the category of topological spaces. Using other algebraic theories gets topological monoids and topological rings. A group in the category of manifolds will be a Lie group. One can go on with these many different examples, e.g., a group in the category of coassociative coalgebras is a Hopf algebra. Continuing the analogy of the cookie cutter, we can use our cookie cutter to make cookies in many different types of dough. Here we concentrate on algebras in the category of sets. 

We are not only interested in models/algebras for an algebraic theory. We are also interested in homomorphisms between models/algebras. Let $(A,+,-,e)$ and $(B,+',-',e')$ be two groups and $f:A\lra B$ be a homomorphism between them. That is, $f$ is a set map that ``respects all the operations'':
\begin{itemize}
\item for all $a_1, a_2 \in A, f(a_1+a_2)=f(a_1)+'f(a_2)$.
\item for all $a_1 \in A, f(-a_1)=-'f(a_1)$.\footnote{For groups, this rule is satisfied if the other two are satisfied} And 
\item for the identity element $e$, $f(e)=e'$. 
\end{itemize}

It will be useful to see these three conditions as commutative diagrams:

$$
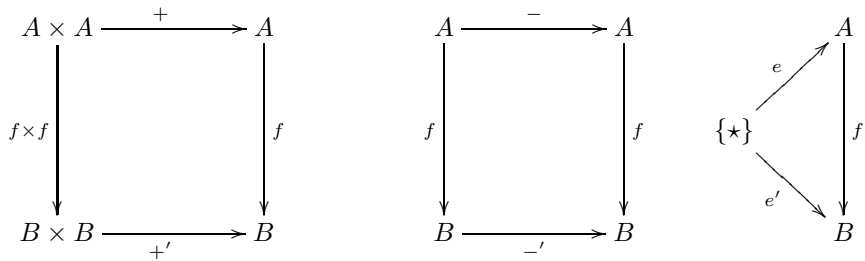
\xymatrix{ 
A \times A \ar[rr]^+ \ar[dd]_{f \times f} && A \ar[dd]^f && A\ar[dd]_f \ar[rr]^{-} &&A\ar[dd]^f&&A\ar[dd]^f\\
&&&&&&&\{ \star \} \ar[ur]^{e}\ar[dr]_{e'}&\\
B \times B \ar[rr]_{+'}&&B && B\ar[rr]_{-'}&&B&& B }$$

If one looks at the two groups as algebras for the algebraic theory of $T_{group}$ in sets, that is, $A:T_{group} \lra \Set$ and $B:T_{group} \lra \Set$, it turns out that a homomorphism from $A$ to $B$ is nothing more than a natural transformation from the functor A to the functor B. The commuting squares shown above is exactly the requirement for a natural transformation. 
For an arbitrary algebraic structure, a homomorphism $f:A \lra B$ is a function between two underlying sets that respects the operations. In terms of algebraic theories, a homomorphism is simply a natural transformation between two functors that are algebras. If an algebraic theory has an $\alpha: n \lra 1$ operation, then a homomorphism $f:(A,\alpha) \lra (B,\alpha')$ 
must `respect' this operation, i.e., 

$$f(\alpha(a_1,a_2, \ldots, a_n))=\alpha'(f(a_1), f(a_2), \ldots, f(a_n)).$$
In terms of commutative diagrams this says that 
$$\xymatrix{ 
A^n \ar[rr]^\alpha \ar[dd]_{f^n} && A \ar[dd]^f \\
&&&&&&\\
B^n \ar[rr]_{\alpha'}&&B }$$ commutes. 

We will need to be slightly more general. If the algebraic theory has an operation from $\alpha:n \lra m$ then a homomorphism must satisfy 
$$\xymatrix{ 
A^n \ar[rr]^\alpha \ar[dd]_{f^n} && A^m \ar[dd]^{f^m} \\
&&&&&&\\
B^n \ar[rr]_{\alpha'}&&B^m. }$$

For an algebraic theory $T$, one can form a category of algebras and homomorphisms between those algebras. We denote this category as $Alg(T,\Set)$. From the categorical point of view this is a subcategory of the functor category $\Set^T$ of functors (and natural transformations) from $T$ to $\Set$. $Alg(T, \Set)$ is the subcategory of functors and natural transformations that respect the product structure of $T$ and $\Set$.

(We should point out that categorical algebra gets far more interesting when we talk about {\em theory morphisms.}These are special types of functors from one algebraic theory $T_1$ to another algebraic theory $T_2$. $F:T_1 \lra T_2$ induces a functor of their categories of algebras $F^*: Alg(T_2, \Set) \lra Alg(T_1, \Set)$. One of the main tools of categorical algebra is that through a powerful tool called a Kan Extension, there is a functor $F_*:Alg(T_1, \Set)\lra Alg(T_2,\Set)$. Much work is done classifying the different types of theory morphisms and the relationships they induce between $F^*$ and $F_*$. It turns out that many theorems in algebra are simple consequences of these ideas. Categorical algebra is one of three ways that mathematicians describe general algebraic structures---the other two being universal algebra and model theory. However, categorical algebra distinguishes itself by being the only one that discusses the relationship between different types of structures. Although this is important for general mathematics, it will not be needed for our discussion. )

For every algebraic theory $T$, there is a forgetful functor from a category of algebras of $T$ to the category of sets, $U:Alg(T, \Set) \lra \Set$. This functor takes an algebra $A:T \lra \Set$ to its underlying set. That is, $U$ is defined as $U(A)=A(1)$ where $A(1)$ is the set that $A$ is defined on. $U$ also takes a homomorphism between two algebraic structures to its underlying set map, $U(f:A \lra B)=f_1:A(1)\lra B(1)$.

We still need some more tools before we can get to our main reconstruction idea. Let $Alg$ be a category which is thought of as a category of algebras of some theory. Let $U:Alg \lra \Set$ be a functor which we will think of as a way of forgetting the algebraic structure. Consider a natural transformation 

$$\alpha: U \Lra U.$$ 

\noindent For every algebra $A$ in $Alg$, there will be a function of sets
$$\alpha_A:U(A) \lra U(A).$$

\noindent This needs to be generalized. Consider a functor $U^n:Alg \lra \Set.$ This functor takes an algebra $A$ to $U^n(A)=(U(A))^n$, the $n$-th product of $U(A)$. A natural transformation 

$$\alpha:U^n \Lra U^m$$

\noindent will associate to every algebra $A$ in $Alg$ a set function

$$\alpha_A:U^n(A) \lra U^m(A).$$

Since $\alpha$ is a natural transformation, it must satisfy the following naturality condition: for every morphism in $Alg$, that is for every homomorphism $f:A \lra B$ in $Alg$, the following square must commute

$$\xymatrix{
U^n(A)\ar[rr]^{\alpha_A}\ar[dd]_{U^n(f)}&& U^m(A)\ar[dd]^{U^m(f)}
\\&&&(\star\star\star)
\\
U^n(B)\ar[rr]_{\alpha_B}&& U^m(B)
}$$

That is basically the setup of categorical algebra. Now comes the magic. We showed how to go from the algebraic theory, $T$ to the category of algebras, $Alg(T, \Set)$ and its forgetful functor $U:Alg(T,\Set) \lra \Set$. Now we go the other way: from a forgetful functor $U:Alg \lra \Set$ we reconstruct the algebraic theory which generates that category of algebras. \footnote{This will not work for all forgetful functors. We need to restrict to functors called ``tractable.'' These are functors $U$ where $Nat(U^n, U^m)$ is a proper set.}

We will construct a theory $T_U$ from $U:Alg \lra \Set$. The objects of $T_U$ are the same as every algebraic theory, the natural numbers. The morphisms (operations) in $T_U$ are defined as follows:

$$T_U(n,m) = Nat(U^n, U^m)$$
where $Nat$ signifies the collection of natural transformations from the functor $U^n$ to $U^m$. It is a major theorem of categorical algebra that this is the algebraic theory we were looking for. If we start off with an algebraic theory $T$ and form the forgetful functor $U:Alg(T, \Set) \lra \Set$, then we can form $T_U$. It is not hard to show that 
$T=T_U$. We have reconstructed the algebraic theory from the category of algebras. 

The philosophy behind this reconstruction is very important for our point. The operations of $T_U$ will correspond to natural transformations $\alpha: U^n \lra U^m$. Let us carefully examine the commutative square $(\star\star\star)$. The usual view is that the homomorphisms $f$ (vertical morphisms) are set functions that respect all the operations $\alpha$ (horizontal morphisms). Now look at the way the operations are defined: the operations are natural transformations, that is, those functions that respect all the homomorphisms (of $Alg$). To reiterate, the usual motto is that 
\centerline{\underline{ Homomorphisms are functions that respect all the operations.}}

\noindent We stress how it works the other way:

\centerline{\underline{Operations are functions that respect all the homomorphisms.}}

Homomorphisms are ways of transforming elements from one set to another.\footnote{Algebraic theories deal with all homeomorphisms, while our paper deals with isomorphisms. The difference is not significant for our needs.} In a sense they are ways of changing the semantics of what we are talking about. An algebra homomorphism must respect the operations of the structure. But for our purposes we can say that the operations are those functions that respect all the possible ways of changing the elements of the sets. Yet another way of saying this is that the operations are those functions that are invariant to any changes that can be made to the elements. Or in our language, operations are functions that respect symmetry of semantics. In a sense, this paper is a vast generalization from universal algebra/category theory to all of mathematics.

Although we highlighted the reconstruction of an algebraic theory from a category of algebras, such reconstruction theories are common in many different parts of modern mathematics. Mathematicians have shown how to 
\begin{itemize} 
\item Reconstruct certain groups from their categories of representations.
\item Reconstruct certain rings from their categories of modules 
\item Reconstruct certain quantum groups (\cite{kassel1995quantum}, \cite{majid2000foundations}, \cite{chari1995guide}) from their monoidal categories of representations, etc.
\end{itemize} Many of these can be seen as instances of the Tannaka Krien Duality Theorem (see \cite{Joyal1991}). As of late, there are also applications of these ideas to quantum field theories. 
While the reconstructions in all these areas are vastly more complicated than what we did above, the main idea that we outlined is at the core of their work. 

All these ideas come from much older ideas that perhaps can be traced back to Felix Klein's Erlangen Program. Namely by looking at the symmetries of a geometric object we can determine properties of that object. Klein was originally only interested in geometric objects, but mathematicians have taken his ideas in many other directions. They look at the automorphism groups and of anything that moves, e.g., groups, models of arithmetic, vector spaces, to algorithms \cite{YanofskyG}. In a sense, some of these ideas can be seen as going back to Galois Theory.

\bibliography{meirbib,universe}
\bibliographystyle{alpha}

\end{document}